\documentclass[reqno,12pt]{amsart}
\usepackage[english]{babel}
\usepackage[dvipsnames]{xcolor}
\usepackage{amssymb,bm}
\usepackage[margin=3cm]{geometry}
\RequirePackage{mathrsfs}

\usepackage[
    colorlinks=true,
    linkcolor=MidnightBlue,
    citecolor=MidnightBlue,
    urlcolor=MidnightBlue
]{hyperref}

\numberwithin{equation}{section}
\renewcommand{\phi}{\varphi}
\renewcommand{\rho}{\varrho}
\renewcommand{\leq}{\leqslant}
\renewcommand{\le}{\leqslant}
\renewcommand{\geq}{\geqslant}
\renewcommand{\ge}{\geqslant}

\newcommand{\PP}{\mathbb{P}}
\newcommand{\QQ}{\mathbb{Q}}
\newcommand{\ZZ}{\mathbb{Z}}
\newcommand{\ZZp}{\mathbb{Z}_{\mathrm{prim}}}
\newcommand{\RR}{\mathbb{R}}
\newcommand{\VV}{\mathbb{V}}
\newcommand{\cA}{\mathcal{A}}
\newcommand{\cB}{\mathcal{B}}
\newcommand{\cC}{\mathcal{C}}
\newcommand{\cR}{\mathcal{R}}
\newcommand{\fd}{\mathfrak{d}}
\newcommand{\x}{\mathbf{x}}
\renewcommand{\c}{\mathbf{c}}
\newcommand{\y}{\mathbf{y}}
\newcommand{\z}{\mathbf{z}}
\renewcommand{\a}{\mathbf{a}}
\renewcommand{\b}{\mathbf{b}}
\renewcommand{\v}{\mathbf{v}}
\renewcommand{\t}{\mathbf{t}}
\newcommand{\X}{\mathbf{X}}
\newcommand{\ee}{\mathbf{e}}
\newcommand{\0}{\mathbf{0}}
\newcommand{\one}{\mathbf{1}}
\newcommand{\norm}[1]{\left\Vert#1\right\Vert}
\newcommand{\loc}{\mathrm{loc}}
\newcommand{\good}{\mathrm{g}}
\DeclareMathOperator{\Span}{Span}
\DeclareMathOperator{\vol}{vol}
\DeclareMathOperator{\rad}{rad}
\DeclareMathOperator{\Br}{Br}
\DeclareMathOperator{\Disc}{Disc}
\newcommand{\Pic}{\operatorname{Pic}}
\newcommand{\Sbar}{\overline S}
\newcommand{\Gal}{\operatorname{Gal}}
\newcommand{\Lzero}{\mathsf{L}_0}

\newtheorem{theorem}{Theorem}[section]
\newtheorem{lemma}[theorem]{Lemma}
\newtheorem{proposition}[theorem]{Proposition}
\newtheorem{corollary}[theorem]{Corollary}
\theoremstyle{remark}

\author{Tim Browning}
\author{Anna Zanoli}

\address{ISTA\\
Am Campus 1\\
3400 Klosterneuburg\\
Austria}
\email{tim.browning@ista.ac.at, anna.zanoli@ista.ac.at}

\address{Department of Mathematics\\ Imperial College London\\
SW7 2AZ\\ 
United Kingdom}
\email{tim.browning@imperial.ac.uk}

\title[Cubic surfaces with a rational point]
      {$\bm{99.99\%}$ of cubic surfaces have a rational point}

\begin{document}
\begin{abstract}
When cubic surfaces in $\PP^3$ over $\QQ$ are ordered by height,
we prove that almost every such surface satisfies the Hasse principle.
This extends work of Browning, Le Boudec and Sawin and completes the proof of the conjecture of Poonen and Voloch for Fano hypersurfaces. We also obtain
an upper bound for the least height of a rational point on almost all everywhere locally soluble cubic surfaces and we give an upper bound for the number of cubic surfaces with a Brauer--Manin
obstruction to the Hasse principle.
Finally, once combined with a calculation 
 of Beneish and Keyes about the density of everywhere locally soluble cubic surfaces, this 
establishes the explicit proportion of cubic surfaces that have a rational point.
\end{abstract}
\maketitle
\thispagestyle{empty}
 \setcounter{tocdepth}{1}
 \tableofcontents

\section{Introduction}

The Hasse principle asks whether the existence of points over every
completion of $\QQ$ guarantees the existence of a rational point.
For quadrics this is the Hasse--Minkowski theorem, whereas cubic
surfaces can fail the Hasse principle. Poonen and Voloch~\cite{poonen}
proposed a statistical formulation for hypersurfaces: in the Fano range,
almost every hypersurface of a fixed degree and dimension should satisfy
the Hasse principle when ordered by the size of its coefficients.
Browning, Le Boudec and Sawin~\cite[Theorem~1.1]{BLS} proved this
prediction in every case except cubic surfaces. Our purpose is to
establish the remaining case.

We first specify the ordering. Listing the twenty monomials of degree
three in four variables lexicographically identifies the space of cubic
surfaces over $\QQ$ with $\VV=\PP^{19}(\QQ)$. For $V\in\VV$, let
$\a_V\in\ZZp^{20}$ be either primitive coefficient vector defining $V\subset \PP^3$,
and let $\norm{\cdot}$ denote the Euclidean norm. Put
\[
 \VV(A)=\{V\in\VV:\norm{\a_V}\leq A\},\qquad
 \rho(A)=
 \frac{\#\{V\in\VV(A):V(\QQ)\ne\emptyset\}}{\#\VV(A)}.
\]
Thus we count equations up to multiplication by a non-zero rational
scalar, rather than isomorphism classes.
The singular members form a set of density zero, so restricting the
statements below to smooth cubic surfaces gives the same limits.

Write $\VV^{\loc}(A)$ for the members of $\VV(A)$ having a point over
$\RR$ and over $\QQ_p$ for every prime $p$. It follows from work of Poonen and
Voloch~\cite[Theorem~3.6]{poonen} that the limit
\[
 \rho^{\loc}
 =\lim_{A\to\infty}\frac{\#\VV^{\loc}(A)}{\#\VV(A)}
\]
exists, is positive, and is a product of local densities.
Indeed, the extension to our Euclidean ordering carries over easily, as explained in 
work of Beneish and Keyes \cite[Remark~2.2]{chris}.
Our main result is the following.

\begin{theorem}\label{thm:main}
We have
\[
 \lim_{A\to\infty}\rho(A)=\rho^{\loc}.
\]
\end{theorem}

In fact, we shall prove in Section \ref{sec:deduction} that the proportion of cubic surfaces of height at most $A$ which
fail the Hasse principle is $O((\log A)^{-1/144})$.
Together with~\cite[Theorem~1.1]{BLS}, this proves the Poonen--Voloch conjecture \cite[Conjecture~2.2(ii)]{poonen} throughout the Fano range, for hypersurfaces ordered by the Euclidean norm of their coefficient vectors, as in~\cite{BLS}. (Poonen and Voloch use the maximum norm; the limit $\rho^{\loc}$ is the same for both orderings by \cite[Remark~2.2]{chris}.)
Beneish and Keyes \cite{chris} give an explicit Euler product for 
$\rho^{\loc}$, and they evaluate it numerically with rigorous error bounds, finding $\rho^{\loc}\approx0.999927$. In fact, Theorem \ref{thm:main} establishes 
\cite[Conjecture~1.2]{chris}, with the following consequence. 

\begin{corollary}
When ordered by height, 
more than $99.99$\% of cubic surfaces over $\QQ$ have a rational point.
\end{corollary}

It was already   observed in \cite{BLS} that 
a positive proportion of cubic surfaces have
a rational point. This  follows from applying  work of Bhargava~\cite[Theorem~2]{Bhargava} 
on plane cubic curves
to a fixed plane section. Furthermore, work of Poonen and
Voloch~\cite[Proposition~3.4]{poonen}  showed that the conclusion of Theorem \ref{thm:main} is true subject to the   hypothesis that the 
Brauer--Manin obstruction explains all failures of the Hasse principle for cubic surfaces. 
This hypothesis was conjectured to be true for all geometrically rational surfaces by 
Colliot-Th\'el\`ene and Sansuc \cite{CTS} in 1979. 
Our next result shows that the Brauer--Manin obstruction  rarely occurs.

\begin{theorem}\label{thm:BM}
We have
\[
\frac{\#\bigl\{V\in\mathbb{V}^{\mathrm{loc}}(A):\ V(\mathbb{A}_{\QQ})^{\Br}=\emptyset\bigr\}}
     {\#\mathbb{V}^{\mathrm{loc}}(A)}
\ll\frac{1}{A}.
\]
\end{theorem}

The restriction to everywhere locally soluble surfaces in the statement of the theorem is natural: if
$V(\mathbb{A}_{\QQ})=\emptyset$ then trivially $V(\mathbb{A}_{\QQ})^{\Br}=\emptyset$,
independently of the Brauer group of $V$.
On the other hand, by work of Poonen and Voloch \cite[Theorem~3.6]{poonen}, the density
$\varrho^{\mathrm{loc}}$ of $\mathbb{V}^{\mathrm{loc}}$ exists and is positive, so
$\#\mathbb{V}^{\mathrm{loc}}(A)\asymp\#\mathbb{V}(A)\asymp A^{20}$. In particular, 
Theorem~\ref{thm:BM} remains true with $\#\mathbb{V}(A)$ in the denominator.

For
$x\in\PP^3(\QQ)$ represented by $\x\in\ZZp^4$, put
$H(x)=\norm{\x}$. This is the anticanonical height on a smooth cubic
surface associated with the Euclidean norm. Define
\[
 \mathfrak M(V)=
 \begin{cases}
 \min_{x\in V(\QQ)}H(x), &\text{ if $V(\QQ)\neq \emptyset$,}
\\
+\infty, &\text{ if $V(\QQ)= \emptyset$.}
\end{cases} 
\]
We obtain the following upper bound for the typical size of the smallest rational point. 

\begin{theorem}\label{thm:height}
Let $\psi:\RR_{>0}\to\RR_{>0}$ satisfy $\psi(u)/u\to\infty$ as
$u\to\infty$. Then
\[
 \lim_{A\to\infty}
 \frac{\#\{V\in\VV(A):
              \mathfrak M(V)\leq\psi(\norm{\a_V})\}}
      {\#\VV(A)}
 =\rho^{\loc}.
\]
\end{theorem}

This is the cubic surface counterpart of~\cite[Theorem~2.2]{BLS}.
In earlier work, 
Le Boudec~\cite[Theorem~1]{LeBoudec} proved that, for any positive
function $\theta(u)=o(u)$, the surfaces satisfying
$\mathfrak M(V)\leq\theta(\norm{\a_V})$ have density zero.
Consequently, for any $\xi(u)\to\infty$, a proportion tending to one
of everywhere locally soluble cubic surfaces satisfies
\[
 \frac{\norm{\a_V}}{\xi(\norm{\a_V})}
 <\mathfrak M(V)\leq\norm{\a_V}\xi(\norm{\a_V}).
\]
In this sense Theorem~\ref{thm:height} gives the optimal scale for the
least height.

\subsection*{Outline of the argument}

We begin by discussing the proof of  Theorems~\ref{thm:main} and \ref{thm:height}, which will eventually be  proved in Section~\ref{sec:deduction}.
Let
$$
 N_V(B)=\#\{x\in V(\QQ):H(x)\leq B\},
$$
for any cubic surface $V\in\VV$.
The method of~\cite{BLS} compares $N_V(B)$ with a localised counting
function $N_V^{\loc}(B)$, defined in Section~\ref{sec:inputs}.
The treatment of the typical size of $N_V^{\loc}(B)$ in~\cite[Proposition~2.4]{BLS} already applies
to cubic surfaces. The main difficulty arises in  the second moment estimate for the
difference $N_V(B)-N_V^{\loc}(B)$.
Indeed, the variance estimate~\cite[Proposition~4.1]{BLS} yields
\begin{equation*} 
\frac1{\# \mathbb{V}_{d,n}(A)} \sum_{V \in \mathbb{V}_{d,n}(A)} \left( N_V(B) - N_V^{\mathrm{loc}}(B) \right)^2 \ll \frac{B}{A},
\end{equation*}
for the  family $\mathbb{V}_{d,n}\subset \PP^{\binom{n+d}{d}-1}(\QQ)$ of degree $d$ Fano hypersurfaces 
 in $\PP^{n}$, provided that  $(d,n) \notin \{(2,2), (3,3) \}$ and  $B/(\log B)^{1/2} \leq A \leq B^2$. 
Here the height used in $N_V(B)$ for a degree $d$ hypersurface in
$\PP^n$ is $\norm{\x}^{n+1-d}$ for a primitive representative $\x$.
(For cubic surfaces this is simply $\norm{\x}$.)
 This cannot extend to the case $(d,n)=(3,3)$ of cubic surfaces in its current form. To see this, fix the line $L$ given by $X_0=X_1=0$.
The condition $L\subset V$ sets four coefficients equal to zero, leaving
$\gg A^{16}$ primitive coefficient vectors of norm bounded by 
$A$. Each corresponding surface has $N_V(B)\gg B^2$. 
On the other
hand, we will see that 
$N_V^{\loc}(B)\ll_\varepsilon B^{1+\varepsilon}$, for any $\varepsilon>0$, 
uniformly for these surfaces.
Thus
\[
 \frac1{\#\VV(A)}\sum_{V\in\VV(A)}
       (N_V(B)-N_V^{\loc}(B))^2
 \gg\left(\frac BA\right)^4,
\]
for all sufficiently large $A,B$. Taking $B= A\sqrt{\log A}$ clearly contradicts the above variance estimate. 
Inspired by recent work of Zanoli \cite{anna} on random quartic surfaces, 
we shall overcome this difficulty by restricting $N_V(B)$ and $N_V^{\loc}(B)$ to \emph{good} points $x=[\x]$, namely those for which the lattice
\[
\Lambda_\x=\{\b\in\ZZ^4:\b\cdot\x=0\}
\]
of linear forms vanishing at $x$ has small largest successive minimum.
This ensures that $x$ lies on no rational line of small height compared with $\|\x\|$, and therefore almost all points of $L$ which were responsible for the large values of $N_V(B)$ in the example above are discarded. 
The argument of \cite{BLS} requires information about the shape of the  lattices $\Lambda_{\nu(\x)}$ and $\Lambda_{\nu(\x)}\cap\Lambda_{\nu(\y)}$, 
for linearly independent $\x,\y\in \ZZp^4$, where 
$\nu:\RR^4\to\RR^{20}$ is the degree three  Veronese map.
Lemma \ref{lem:products} bounds their successive minima in terms of the largest successive minimum of 
 $\Lambda_\x$ and $\Lambda_\y$, and these bounds are small when $x$ and $y$ are both good.
 Discarding points can only decrease $N_V(B)$, which is harmless since we only need a lower bound for it. Lemma~\ref{lem:discard} shows that the non-good points have sufficiently small total weight
to make their contribution to the local count negligible, uniformly
away from coefficient vectors of small norm. In Section~\ref{sec:moments} we then show that the three truncated moments have the same main term, which cancels in the variance. All the remaining tools and estimates are supplied by~\cite{BLS}. Finally, Section~\ref{sec:deduction} will wrap up the proof of 
Theorems~\ref{thm:main} and \ref{thm:height}.

Section \ref{s:BMO} is concerned with the proof of Theorem \ref{thm:BM}.
While this paper was being completed, O'Dorney
\cite{ODorney} posted a preprint on the Galois group of the $27$ lines of a
random cubic surface. 
Our Proposition~\ref{prop:Br} below can be deduced from
\cite[Proposition~3.1]{ODorney},
together with results of Elsenhans and Jahnel 
\cite{EJ-BM, EJ-order3} that classify which configurations of lines give rise to non-constant Brauer groups for smooth cubic surfaces. 
Since the proof of \cite[Proposition~3.1]{ODorney} is rather 
brief, we have chosen to give a detailed treatment 
in Section \ref{s:BMO}, in which we construct the relevant covers explicitly
and verify that they are integral.

\subsection*{Notation}

Throughout, implied constants are absolute unless indicated otherwise,
and $A$ and $B$ are sufficiently large. We write $\cB_j(u)$ for the
closed Euclidean ball of radius $u$ in $\RR^j$, $\one_E$ for the
indicator of a condition $E$, and $\rad(q)$ for the product of the
distinct prime divisors of a positive integer $q$.
We denote the volume of $\cB_j(1)$ by $C_j$. The elementary asymptotic
\[
 \#\VV(A)\sim\frac{C_{20}}{2\zeta(20)}A^{20}
\]
will be used throughout; it follows from a simple geometry of numbers argument (as
explained by Le Boudec \cite[Lemma~3]{LeBoudec}, for example).


\subsection*{AI statement}
The arguments in Section \ref{s:BMO} were developed with the assistance of ChatGPT-$6$ and Claude Opus $5.5$; all details have been checked by the authors, who take responsibility for them.

\section{Bounds for successive minima}

For a non-zero integral vector $\v\in\ZZ^m$, put
\[
 \Lambda_{\v}=\{\b\in\ZZ^m:\b\cdot\v=0\}.
\]
We write $\lambda_i(\Lambda)$ for the successive minima of a lattice $\Lambda$
with respect to the Euclidean unit ball. Let
$\nu=\nu_{3,3}:\RR^4\to\RR^{20}$ be the Veronese map listing the
degree three monomials in lexicographical order. 
For $\x\in\ZZ^4_{\text{prim}}$, the lattice $\Lambda_{\nu(\x)}\subset \ZZ^{20}$ is a  
primitive lattice of rank $19$.
For linearly independent
primitive vectors $\x,\y$, put
\begin{equation}\label{eq:Gamma}
 \Gamma_{\x,\y}=\Lambda_{\nu(\x)}\cap\Lambda_{\nu(\y)}\subset \ZZ^{20}.
\end{equation}
This is a primitive lattice of rank $18$.
Building on \cite[Lemma~3.15]{BLS} and 
\cite[Lemma~2.4]{anna}, 
the following result provides upper bounds for the successive minima of these lattices.

\begin{lemma}\label{lem:products}
For  $\x,\y\in\ZZp^4$ such that $\x\neq\pm \y$, we have
\begin{equation}\label{eq:product-minimum}
 \lambda_{18}(\Gamma_{\x,\y})
 \ll \lambda_3(\Lambda_{\x})\lambda_3(\Lambda_{\y}).
\end{equation}
Moreover,
\begin{equation}\label{eq:single-minimum}
 \lambda_{19}(\Lambda_{\nu(\x)})\leq
 \lambda_3(\Lambda_{\x}).
\end{equation}
\end{lemma}

\begin{proof}
Choose three linearly independent integral linear forms
$L_1,L_2,L_3$ vanishing at $\x$, whose coefficient vectors have
norm at most $\lambda_3(\Lambda_{\x})$. Choose
$M_1,M_2,M_3$ similarly for $\y$.
For $1\leq i,j\leq3$ and $0\leq k\leq3$, 
the cubics
\[
 L_i(\X)M_j(\X)X_k
\]
have integral coefficients, vanish at both points, and have
coefficient norms
\[
 \ll \lambda_3(\Lambda_{\x})\lambda_3(\Lambda_{\y}).
\]
It therefore suffices to show that their real span has
dimension eighteen.

Since the $L_i$ and $M_j$ span the spaces of linear forms
vanishing at $\x$ and $\y$, respectively, this dimension
can be computed after an invertible real change of coordinates.
Take coordinates in which $\x=\ee_0$ and $\y=\ee_1$.
The two spaces of linear forms are then
\[
 \langle X_1,X_2,X_3\rangle
 \quad\text{and}\quad
 \langle X_0,X_2,X_3\rangle.
\]
Their products with arbitrary linear forms span exactly
the cubic monomials other than $X_0^3$ and $X_1^3$. Indeed, 
each of the other eighteen monomials can be written as
a product of one variable from each displayed list and
one further variable.
Thus the original integral products span an
eighteen-dimensional space. Selecting eighteen independent
ones proves~\eqref{eq:product-minimum}.
(Note that the change of coordinates was used only for this dimension
calculation, so it does not affect the norm bound.)

For~\eqref{eq:single-minimum}, consider the integral cubics
$L_i(\X)Q$, where $Q$ runs over the quadratic monomials.
Multiplication by a monomial preserves the coefficient norm,
so each has norm at most $\lambda_3(\Lambda_{\x})$.
In coordinates with $\x=\ee_0$, these products span all
cubic monomials except $X_0^3$. Their span therefore has
dimension nineteen. Selecting nineteen independent products
in the original coordinates proves the second bound.
\end{proof}

For $\x\in\ZZ^4\setminus \{\0\}$, let $\fd_2(\x)$ be the
minimum determinant of a rank two sublattice of $\ZZ^4$ containing
$\x$: 
\[
\mathfrak{d}_2(\x)
=
\min\left\{
\det(\mathcal{L}):
\mathcal{L}\subset\mathbb{Z}^4,\ 
\operatorname{rank}(\mathcal{L})=2,\ 
\x\in \mathcal{L}
\right\}.
\] Equivalently, $\mathfrak{d}_2(\x)$ is the least height of a rational line
through $x$, where the height of a rational line $\ell=\mathbb{P}(\Pi)$, for
a plane $\Pi\subset\RR^4$, is $\det(\Pi\cap\ZZ^4)$, the Euclidean norm of its
primitive Pl\"ucker vector. Put
\[
 \mu(\x)=\frac{3\norm{\x}}{\fd_2(\x)}.
\]
Taking $d=1$ and $n=3$ in~\cite[Lemma~3.15]{BLS} yields $\lambda_3(\Lambda_{\x})\leq\mu(\x)$.
Adopting the terminology from \cite[Section~3]{anna}, we introduce the  collection of {\em good} points to be the  set
 \begin{equation}\label{eq:good-set}
  \Xi^{\good}(B,T)=\{\x\in\ZZp^4:
               \norm{\x}\leq B,\ \mu(\x)\leq T\},
\end{equation}
for  $T\geq 1$. Note that $\mathfrak{d}_2(\x)\leq\|\x\|$, since $\ZZ\x+\ZZ\mathbf{e}_i$ has determinant at most $\|\x\|$ for any standard basis vector $\mathbf{e}_i\notin\RR\x$. Hence $\mu(\x)\geq3$, and $\Xi^{\good}(B,T)$ is empty when $T<3$.
This set  is invariant under $\x\mapsto-\x$. Uniformly for $\x,\y\in\Xi^{\good}(B,T)$ with 
$\x\ne\pm\y$,
Lemma~\ref{lem:products}
gives
\begin{equation}\label{eq:good-minima}
 \lambda_{19}(\Lambda_{\nu(\x)})\leq T,\qquad
 \lambda_{18}(\Gamma_{\x,\y})\ll T^2.
\end{equation}
The next result shows that there is a negligible contribution from the average of a key quantity over non-good points. 

\begin{lemma}\label{lem:discard}
For $B,T\geq1$, we have
\[
 \sum_{\substack{\x\in\ZZp^4\setminus\Xi^{\good}(B,T)\\
                         \norm{\x}\leq B}}
 \frac1{\norm{\nu(\x)}}\ll \frac{B^2}{T^3}.
\]
\end{lemma}

\begin{proof}
For $X,\Delta\geq1$, it follows from \cite[Lemma~2.6]{anna} that 
\begin{equation}\label{eq:small-d2}
 \#\{\x\in\ZZ^4:0<\norm{\x}\leq X,
                       \fd_2(\x)\leq\Delta\}
 \ll X^2\Delta^3,
\end{equation}
refining \cite[Lemma~3.20]{BLS}.
In a dyadic interval $X/2<\norm{\x}\leq X$, a discarded point
satisfies $\fd_2(\x)<3X/T$. Since an integral rank two lattice has
determinant at least one, there are no such points if $3X<T$.
Otherwise~\eqref{eq:small-d2} bounds their number by $O(X^5T^{-3})$.
Using $\norm{\nu(\x)}\asymp\norm{\x}^3$ and summing over dyadic
$X\ll B$, we obtain
\[
\sum_{\substack{\x\in\ZZp^4\setminus\Xi^{\good}(B,T)\\
                         \norm{\x}\leq B}}
\frac1{\norm{\nu(\x)}}
 \ll\frac1{T^3}\sum_{\substack{X\ll B\\
 \text{dyadic}}
 }X^2
 \ll\frac{B^2}{T^3},
\]
as required. 
\end{proof}

\section{The local count and auxiliary estimates}
\label{sec:inputs}

As in~\cite[Section~2]{BLS}, put
\[
 \alpha=\log B,\qquad w=\frac{\log B}{\log\log B},\qquad
 W=\prod_{p\leq w}p^{\lceil\log w/\log p\rceil+1}.
\]
It follows from \cite[Equation (2.12)]{BLS} that $W\ll B^{4/\log\log B}$. Hence any fixed power of
$\alpha W$ is $O_\varepsilon(B^{\varepsilon})$. 
For $\v\in\RR^{20}$ and an integer
$q\geq1$, let
\begin{align*}
 \cC_{\v}^{(\alpha)}
 &=\left\{\t\in \RR^{20}:|\v\cdot\t|
           \leq\frac{\norm{\v}\norm{\t}}{2\alpha}\right\},\\
 \Lambda_{\v}^{(q)}
 &=\{\b\in\ZZ^{20}:\v\cdot\b\equiv0\bmod q\},
\end{align*}
where the second definition is used only for integral $\v$. Write 
\[
\Xi(B)=\{\x\in\ZZp^4:\norm{\x}\leq B\},
\]
so that $\Xi^\good(B,T)\subset \Xi(B)$ in \eqref{eq:good-set}, for any $T\geq 1$. Define the localised counting function as 
$$N_V^{\mathrm{loc}}(B)
=
\frac{1}{2}\cdot\frac{\alpha W}{\lVert \a_V\rVert}
\sum_{\substack{\x\in\Xi(B)\\
\a_V\in \Lambda_{\nu(\x)}^{(W)}
\cap \mathcal{C}_{\nu(\x)}^{(\alpha)}}}
\frac{1}{\lVert \nu(\x)\rVert}.$$

For every non-zero $\a\in\ZZ^{20}$ define
\begin{align}
 N_{\a}^{\good}(B,T)
 &=\frac12\sum_{\x\in\Xi^{\good}(B,T)}
                    \one_{\a\cdot\nu(\x)=0},\label{eq:Ng}\\
 L_{\a}^{\good}(B,T)
 &=
 \frac{1}{2}\cdot
 \frac{\alpha W}{\norm{\a}}
   \sum_{\x\in\Xi^{\good}(B,T)}\frac{
      \one_{\a\in\Lambda_{\nu(\x)}^{(W)}
                   \cap\cC_{\nu(\x)}^{(\alpha)}}}
                   {\norm{\nu(\x)}}.\label{eq:Lg}
\end{align}
Denote by $L_{\a}(B)$ 
the analogous local count with $\Xi(B)$ in place of $\Xi^{\good}(B,T)$, so that if $\a$ is either of the two primitive coefficient vectors
for $V$, then $L_{\a}(B)=N_V^{\loc}(B)$.
Lemma~\ref{lem:discard} immediately yields the pointwise estimate
\begin{equation}\label{eq:local-loss}
 0\leq L_{\a}(B)-L_{\a}^{\good}(B,T)
 \ll \frac{\alpha W B^2}{\norm{\a}T^3}.
\end{equation}
We shall ultimately work with $A,B$ in the range
\begin{equation}\label{eq:range}
 A\leq B\leq A(\log A)^{1/2}
\end{equation}
and we shall choose $T=A^{5/12}$.

We record some more  inputs from~\cite{BLS}.
Given linearly independent vectors $\c,\mathbf d\in \ZZ^N$,
let 
$G(\c,\mathbf d)=\gcd_{1\leq i<j\leq N}(c_i d_j-c_j d_i)$.
Let $\x,\y\in\Xi(B)$ with $\x\ne\pm\y$. It follows from
\cite[Lemma~3.12]{BLS} that $G(\nu(\x),\nu(\y))=G(\x,\y)$, and we set
\begin{equation}\label{eq:notation}
 \v=\nu(\x),\quad\z=\nu(\y),\quad
 G=\gcd_{0\leq i,j\leq 3}(x_i y_j-x_j y_i),\quad
 \Delta=\frac{\norm{\v}\norm{\z}}
                  {\norm{\v\wedge\z}}.
\end{equation}
The vectors $\v$ and $\z$ are primitive, and they are linearly independent since $\x\neq\pm\y$. Moreover, $G>0$ because $\x$ and $\y$ are linearly independent, and $\Delta\geq1$ because $\|\v\wedge\z\|\leq\|\v\|\|\z\|$.
Recalling the definition \eqref{eq:Gamma} of $\Gamma_{\x,\y}$, 
the determinant  identities
in~\cite[Lemmas~3.9 and 3.10]{BLS} give
\begin{align}
 \det\Gamma_{\x,\y}&=\frac{\norm{\v\wedge\z}}{G},\notag \\
 \det(\Lambda_{\v}\cap\Lambda_{\z}^{(W)})
   &=\frac{W\norm{\v}}{\gcd(G,W)},\label{eq:Dmix}\\
 \det(\Lambda_{\v}^{(W)}\cap\Lambda_{\z}^{(W)})
   &=\frac{W^2}{\gcd(G,W)}.\label{eq:Dloc}
\end{align}
In particular, we have 
\begin{equation}\label{eq:D-inverse}
\frac{1}{\det\Gamma_{\x,\y}}
=\frac{G\Delta}{\norm{\v}\norm{\z}}.
\end{equation}

Writing $e_{\x,\y}=\min\{1,\Delta^2/\alpha^2\}$, we claim that 
\begin{align}
 \vol_{19}(\v^{\perp}\cap\cB_{20}(1)
                         \cap\cC_{\z}^{(\alpha)})
 &=\frac{18}{19}C_{18}\frac{\Delta}{\alpha}
                      (1+O(e_{\x,\y})),\label{eq:volume-mix}\\
 \vol_{20}(\cB_{20}(1)\cap\cC_{\v}^{(\alpha)}
                         \cap\cC_{\z}^{(\alpha)})
 &=\frac{18}{20}C_{18}\frac{\Delta}{\alpha^2}
                      (1+O(e_{\x,\y})).\label{eq:volume-loc}
\end{align}
We can obtain these estimates from~\cite[Lemmas~4.2 and~4.3]{BLS}, by taking  $N=20$ and putting
\[
 \v'=\frac{2\alpha\v}{\norm{\v}},
 \qquad
 \z'=\frac{2\alpha\z}{\norm{\z}}.
\]
Indeed, $\t\in\cC_{\v}^{(\alpha)}$ is equivalent to
$|\v'\cdot\t|\leq\norm{\t}$, and similarly for $\z$.
Thus the two volumes above are respectively
$\mathcal I(\v',\z')$ and $\mathcal J(\v',\z')$
in the notation of~\cite[\S4.2]{BLS}. Since $\norm{\v'}=\norm{\z'}=2\alpha$ and 
\[
 \norm{\v'}^2\norm{\z'}^2-(\v'\cdot\z')^2=\norm{\v'\wedge\z'}^2=
 \frac{16\alpha^4}{\norm{\v}^2\norm{\z}^2}\norm{\v\wedge\z}^2
  =\frac{16\alpha^4}{\Delta^2},
\]
by \eqref{eq:notation}, 
the stated main terms and errors follow.

The estimates for inverse determinants in~\cite[Lemmas~4.5 and~4.6]{BLS}
include the case $(d,n)=(3,3)$ and give
\begin{align}
 \sum_{\substack{\x,\y\in\Xi(B)\\\x\ne\pm\y}}
 \frac{1}{\det\Gamma_{\x,\y}}
 &\ll B^2,\label{eq:E-bound}\\
 \sum_{\substack{\x,\y\in\Xi(B)\\\x\ne\pm\y}}
 \frac{e_{\x,\y}+\one_{G\nmid (W/\rad(W))}}{
 \det\Gamma_{\x,\y}}
 &\ll\frac{B^2}{(\log B)^{1/2}}.\label{eq:F-bound}
\end{align}
The  summands in these  bounds are non-negative, so the same upper
bounds hold after restricting to any subset of pairs. Notice also that
$\one_{G\nmid W}\leq\one_{G\nmid (W/\rad(W))}$.

Next,  it follows from~\cite[Proposition~2.4]{BLS} that 
\begin{equation}\label{eq:local-small}
 \frac{\#\{V\in\VV^{\loc}(A):
          N_V^{\loc}(A\phi(A))\leq\phi(A)^{2/3}\}}
      {\#\VV^{\loc}(A)}
 \ll\phi(A)^{-1/72}
\end{equation}
for $1<\phi(A)\leq A$. 

\medskip

Finally, 
we state explicitly a weighted lattice point counting  estimate 
that will be used in the mixed and local moments. If $\Lambda\subset \ZZ^{20}$ is an
integral lattice of rank $R$, we have the standard bound 
\begin{equation}\label{eq:crude-lattice}
 \#(\Lambda\cap\cB_{20}(u))\ll u^R,
\end{equation}
for any $u\geq 1$, in which the implied constant is independent of $\Lambda$.

\begin{lemma}\label{lem:weighted}
Let $\Lambda\subset\ZZ^{20}$ be a lattice of rank $R$, and let
$\cR(A)$ be the intersection of $\cB_{20}(A)$ with at
most two cones of the form $\cC_{\v}^{(\alpha)}$.
Put
$v_{\cR}=\vol_R(\Span_{\RR}\Lambda\cap\cR(1))$.
Suppose $1\leq Y\leq A$ and $\lambda_R(\Lambda)\leq Y$.
If $k\in\{1,2\}$ and $R-k>1$, then
\[
 \sum_{\a\in(\Lambda\setminus\{\0\})\cap\cR(A)}
       \frac1{\norm{\a}^k}
 =\frac{R}{R-k}\frac{v_{\cR}A^{R-k}}{\det\Lambda}
 +O\left(\frac{YA^{R-k-1}}{\det\Lambda}+Y^{R-k}\right).
\]
\end{lemma}

\begin{proof}
For $Y\le u\le A$, apply~\cite[Lemma~3.4]{BLS} to the ball
and~\cite[Lemma~3.5]{BLS} when one or two cone conditions are
present. Since $\det\Lambda\ll Y^R\le Yu^{R-1}$, the error absorbs
the contribution of the origin. Thus the number of non-zero lattice
points in $\cR(u)$ is
\[
 \frac{v_{\cR}u^R}{\det\Lambda}
   +O\left(\frac{Yu^{R-1}}{\det\Lambda}\right).
\]
Apply partial summation on $[Y,A]$. The leading term is
$\frac R{R-k}\frac{v_{\cR}A^{R-k}}{\det\Lambda}$, up to
$O(Y^{R-k})$, since $\det\Lambda\geq1$. The integrated error is
$O(YA^{R-k-1}/\det\Lambda)$, since $R-k>1$.
The part of the weighted sum with $\norm{\a}\leq Y$ is
$O(Y^{R-k})$ by~\eqref{eq:crude-lattice} and dyadic summation.
\end{proof}

\section{The truncated moment calculation}\label{sec:moments}

Throughout this section we shall assume that $A,B$ are in the range~\eqref{eq:range}, and 
we let $T\in[1,\sqrt{A}]$ be a parameter. 
Set
\[
 \cA(A)=\{\a\in\ZZ^{20}:0<\norm{\a}\leq A\}
 \]
 and 
 \[
 \Omega^{\good}(B,T)=\{(\x,\y)\in\Xi^{\good}(B,T)^2:
                                  \x\ne\pm\y\},
\]
in the notation of \eqref{eq:good-set}.
Write
\[
 E^{\good}(B,T)=\sum_{(\x,\y)\in\Omega^{\good}(B,T)} \frac{1}{\det\Gamma_{\x,\y}}
 \]
and put
\begin{equation}\label{eq:c}
c=\frac{C_{18}}4.
\end{equation}
It follows from~\eqref{eq:E-bound} that $E^{\good}(B,T)\ll B^2$.
We also use the elementary estimates
\begin{equation}\label{eq:elementary-sums}
 \#\Xi(B)\ll B^4,\qquad
 \sum_{\x\in\Xi(B)}\frac{1}{\norm{\nu(\x)}}\ll B,\qquad
 \sum_{\x\in\Xi(B)}\frac{1}{\norm{\nu(\x)}^{2}}\ll1.
\end{equation}

Recall the notation \eqref{eq:Ng} and \eqref{eq:Lg} for 
$N_{\a}^{\good}(B,T)$ and $L_{\a}^{\good}(B,T)$, respectively.
Let $D^{\good}$, $D^{\mathrm{mix},\good}$ and
$D^{\loc,\good}$ be the off-diagonal parts of
\[
 \sum_{\a\in\cA(A)}(N_{\a}^{\good}(B,T))^2,\qquad
 \sum_{\a\in\cA(A)}N_{\a}^{\good}(B,T)L_{\a}^{\good}(B,T),\qquad
 \sum_{\a\in\cA(A)}(L_{\a}^{\good}(B,T))^2,
\]
respectively.
Here off-diagonal means $\x\ne\pm\y$. We sum over all non-zero
coefficient vectors of norm at most $A$. This avoids M\"obius inversion
in the moment calculation and suffices for our application, since the
square of $N_{\a}^{\good}-L_{\a}^{\good}$ is non-negative and each
surface has two primitive representatives.

\subsection{The global moment}

The lattice estimate~\cite[Lemma~3.4]{BLS}, together with
\eqref{eq:good-minima}, gives
\begin{equation}\label{eq:oscar}
 \#\left((\Gamma_{\x,\y}\setminus\{\0\})
                     \cap\cB_{20}(A)\right)
 =\frac{C_{18}A^{18}}{\det\Gamma_{\x,\y}}
  +O\left(\frac{A^{17}T^2}{\det\Gamma_{\x,\y}}\right)
\end{equation}
uniformly on $\Omega^{\good}(B,T)$. Here we have used the fact that
$\lambda_{18}(\Gamma_{\x,\y})\ll T^2\ll A$. Moreover, the
contribution of the origin is absorbed into the error term, since
Minkowski's second theorem gives
$\det\Gamma_{\x,\y}\leq\lambda_{18}(\Gamma_{\x,\y})^{18}\ll T^{36}\leq A^{17}T^2$.
Expanding the square in~\eqref{eq:Ng} and swapping the order of summation, we obtain $$D^{\good}=\frac14\sum_{(\x,\y)\in\Omega^{\good}(B,T)}\#\big((\Gamma_{\x,\y}\setminus\{\0\})\cap\cB_{20}(A)\big).$$
Summing \eqref{eq:oscar} over $\Omega^{\good}(B,T)$, and using the estimate
$E^{\good}(B,T)\ll B^2$, we get
\begin{equation}\label{eq:global-moment}
 D^{\good}=cA^{18}E^{\good}(B,T)
       +O(A^{17}T^2B^2),
\end{equation}
where $c$ is given by \eqref{eq:c}.

\subsection{The mixed moment}

Fix $(\x,\y)\in\Omega^{\good}(B,T)$ and use the notation
$\v,\z,G,\Delta$ from \eqref{eq:notation}.
Put 
\[
 \Gamma^{\mathrm{mix}}=\Lambda_{\v}
                         \cap\Lambda_{\z}^{(W)}.
                         \]
This lattice has rank $19$, since $W\Lambda_{\v}\subseteq\Gamma^{\mathrm{mix}}\subseteq\Lambda_{\v}$ where both lattices have rank $19$. Multiplying $19$ independent vectors of
$\Lambda_{\v}$ of norm at most $T$ by $W$, it follows from \eqref{eq:good-minima} that 
\[
 \lambda_{19}(\Gamma^{\mathrm{mix}})\leq WT.
\]
Moreover,  \eqref{eq:Dmix} yields
                         \[
                          \det \Gamma^{\mathrm{mix}}
=\frac{W\norm{\v}}{\gcd(G,W)}.
\]
By~\eqref{eq:Ng} and~\eqref{eq:Lg}, expanding  the product and swapping the order of summation gives
\begin{equation}\label{eq:Dmix-expand}
 D^{\mathrm{mix},\good}=\sum_{(\x,\y)\in\Omega^{\good}(B,T)}\frac{\alpha W}{4\|\z\|}
 \sum_{\substack{\a\in\mathcal{A}(A)\cap\Gamma^{\mathrm{mix}}\\ \a\in\mathcal{C}^{(\alpha)}_{\z}}}\frac{1}{\|\a\|},
\end{equation}
where $\v$, $\z$ and $\Gamma^{\mathrm{mix}}$ depend on $(\x,\y)$ as above. Applying Lemma~\ref{lem:weighted} with $(R,k,Y)=(19,1,WT)$, and combining it with
\eqref{eq:volume-mix}, gives
\begin{align}
 \sum_{\substack{\a\in\cA(A)\cap\Gamma^{\mathrm{mix}}\\
                        \a\in\cC_{\z}^{(\alpha)}}}
 \frac1{\norm{\a}}
 &=\frac{C_{18}A^{18}}{\alpha \det \Gamma^{\mathrm{mix}}}
                      \Delta(1+O(e_{\x,\y}))+O\left(\frac{WT A^{17}}{
                      \det \Gamma^{\mathrm{mix}}}
                         +(WT)^{18}\right).\label{eq:mixed-pair}
\end{align}
Inserting this into~\eqref{eq:Dmix-expand}, the main term of the summand for $(\x,\y)$ is
\[
 cA^{18}\frac{\gcd(G,W)\Delta}{\norm{\v}\norm{\z}}(1+O(e_{\x,\y}))
       =cA^{18}\frac{\gcd(G,W)}{G}\frac{1}{\det \Gamma_{\x,\y}}(1+O(e_{\x,\y})),
\]
by \eqref{eq:D-inverse},
where $c$ is given by 
\eqref{eq:c}.
Its difference from $cA^{18}/\det \Gamma_{\x,\y}$ is
\[
\ll\frac{A^{18}}{\det \Gamma_{\x,\y}}
                (e_{\x,\y}+\one_{G\nmid W}),
\]
because $0<\gcd(G,W)/G\leq1$, with 
$\gcd(G,W)/G=1$ if and only if $G\mid W$.
The first error term in~\eqref{eq:mixed-pair}, after multiplication by $\frac{\alpha W}{4\|\z\|}$,
is 
\[
\ll \frac{\alpha WT A^{17}}{\det \Gamma_{\x,\y}},
\]
using $\gcd(G,W)\leq G$ and $\Delta\geq1$. 
Once summed over $\x,\y$, the  remaining error term is
\[
\ll \alpha W(WT)^{18}\sum_{\x,\y\in\Xi(B)}\norm{\nu(\y)}^{-1}
 \ll\alpha W(WT)^{18}B^5, 
\]
by~\eqref{eq:elementary-sums}.

It now follows from~\eqref{eq:E-bound} and~\eqref{eq:F-bound} that
\begin{align}
 D^{\mathrm{mix},\good}
 &=cA^{18}E^{\good}(B,T)
 +O\left(\frac{A^{18}B^2}{(\log B)^{1/2}}
          +\alpha WT A^{17}B^2 +\alpha W(WT)^{18}B^5\right).
 \label{eq:mixed-moment}
\end{align}

\subsection{The local moment}

The lattice
\[
 \Gamma^{\loc}=\Lambda_{\v}^{(W)}
                    \cap\Lambda_{\z}^{(W)}
\]
has rank $20$, 
and 
determinant $\det  \Gamma^{\loc}=W^2/\gcd(G,W)$, 
by \eqref{eq:Dloc}. Moreover, 
its largest successive minimum is at most
 $W$, since it contains $W\ZZ^{20}$. Similarly to above, we have
\[
 D^{\loc,\good}=\sum_{(\x,\y)\in\Omega^{\good}(B,T)}\frac{\alpha^2W^2}{4\|\v\|\|\z\|}
 \sum_{\substack{\a\in\mathcal{A}(A)\cap\Gamma^{\loc}\\ \a\in\mathcal{C}^{(\alpha)}_{\v}\cap\mathcal{C}^{(\alpha)}_{\z}}}\frac{1}{\|\a\|^2}.
\]
Lemma~\ref{lem:weighted}, with $(R,k,Y)=(20,2,W)$, and
\eqref{eq:volume-loc} gives
\[
 \sum_{\substack{\a\in\cA(A)\cap\Gamma^{\loc}\\
          \a\in\cC_{\v}^{(\alpha)}
                    \cap\cC_{\z}^{(\alpha)}}}
 \frac1{\norm{\a}^2}
 =\frac{C_{18}A^{18}}{\alpha^2
\det \Gamma^{\loc}
}
                         \Delta(1+O(e_{\x,\y})) 
                         +O\left(\frac{WA^{17}}{\det \Gamma^{\loc}}+W^{18}\right).
\]
Multiplying by $\alpha^2W^2/(4\norm{\v}\norm{\z})$,
the main term is again 
\[
cA^{18}\frac{\gcd(G,W)}{G}\frac{1}{\det \Gamma_{\x,\y}}(1+O(e_{\x,\y})).
\] 
As in the mixed moment treatment, we can approximate this by 
$cA^{18}/\det \Gamma_{\x,\y}$ with error 
$O(A^{18}B^2/(\log B)^{1/2})$, once summed over the remaining $\x,\y$ using \eqref{eq:F-bound}.
 The other errors sum to
\[
\ll \alpha^2W^2A^{17}B^2+\alpha^2W^{20}B^2,
\]
on taking $\det \Gamma^{\loc}\geq W$ in the treatment of the first term. 
Thus
\begin{equation}\label{eq:local-moment}
 D^{\loc,\good}=cA^{18}E^{\good}(B,T)
 +O\left(\frac{A^{18}B^2}{(\log B)^{1/2}}
          +\alpha^2W^2A^{17}B^2+\alpha^2W^{20}B^2\right).
\end{equation}

\subsection{Diagonal terms and the variance}

Because the point set is invariant under sign, the diagonal part of
$(N_{\a}^{\good}(B,T))^2$ is exactly $N_{\a}^{\good}(B,T)$.
By~\eqref{eq:good-minima}
and \cite[Lemma~3.6]{BLS}, 
 a lattice count for
$\Lambda_{\nu(\x)}$, whose determinant is $\norm{\nu(\x)}$, gives
\begin{equation}\label{eq:first-moment}
 \sum_{\a\in\cA(A)}N_{\a}^{\good}(B,T)
 \ll A^{19}\sum_{\x\in\Xi^{\good}(B,T)}\frac{1}{\norm{\nu(\x)}}
 \ll A^{19}B,
\end{equation}
by \eqref{eq:elementary-sums}.
The diagonal part $\Delta_{\a}^{\loc,\good}$ of
$(L_{\a}^{\good}(B,T))^2$ satisfies
\[
 0\leq\Delta_{\a}^{\loc,\good}
 \leq\frac{\alpha^2W^2}{2\norm{\a}^2}
               \sum_{\x\in\Xi(B)}\norm{\nu(\x)}^{-2}
 \ll\frac{\alpha^2W^2}{\norm{\a}^2}.
\]
Hence
\begin{equation}\label{eq:local-diagonal}
 \sum_{\a\in\cA(A)}\Delta_{\a}^{\loc,\good}
 \ll\alpha^2W^2 A^{18}.
\end{equation}
The mixed diagonal is non-negative and occurs with coefficient $-2$
when the square is expanded;  it can therefore be omitted in an upper
bound. In this way, we obtain
\begin{align*}
 \sum_{\a\in\cA(A)}(N_{\a}^{\good}(B,T)-L_{\a}^{\good}(B,T))^2
 &\leq D^{\good}-2D^{\mathrm{mix},\good}+D^{\loc,\good}
 +\sum_{\a\in\cA(A)}
                  (N_{\a}^{\good}(B,T)+\Delta_{\a}^{\loc,\good}).
\end{align*}
The three main terms cancel. Indeed, on combining
\eqref{eq:global-moment}, \eqref{eq:mixed-moment},
\eqref{eq:local-moment}, \eqref{eq:first-moment} and
\eqref{eq:local-diagonal}, we deduce that 
\begin{align*}
 \frac1{A^{20}}\sum_{\a\in\cA(A)}
       (N_{\a}^{\good}(B,T)-L_{\a}^{\good}(B,T))^2
 \ll {}&\frac BA+\frac{B^2}{A^2(\log B)^{1/2}}
       +\frac{(T^2+\alpha WT+\alpha^2W^2)B^2}{A^3}\\
 &+\frac{B^8+\alpha W(WT)^{18}B^5
                    +\alpha^2W^{20}B^2}{A^{20}}
       +\frac{\alpha^2W^2}{A^2},
\end{align*}
for any parameter $T\in[1,\sqrt{A}]$.

We shall take $T=A^{5/12}$ and we recall that 
$(\alpha W)^{20}=O_\varepsilon(B^{\varepsilon})$, for any $\varepsilon>0$.
Assuming that $A,B$ are in the range~\eqref{eq:range}, it therefore follows that 
\[
 \frac{T^2B^2}{A^3}+
 \frac{\alpha WT B^2}{A^3}+
 \frac{\alpha^2W^2B^2}{A^3}+
 \frac{\alpha W(WT)^{18}B^5}{A^{20}}\ll_\varepsilon A^{-1/6+\varepsilon}.
\]
Taking $\varepsilon=1/12$, 
we have proved the following bound.

\begin{proposition}\label{prop:variance}
Uniformly for $A,B$ in the range $A\leq B\leq A(\log A)^{1/2}$, we have
\[
 \frac1{A^{20}}\sum_{\a\in\cA(A)}
   (N_{\a}^{\good}(B,A^{5/12})-L_{\a}^{\good}(B,A^{5/12}))^2
 \ll\frac BA+\frac{B^2}{A^2(\log B)^{1/2}}+A^{-1/12}.
\]
\end{proposition}

\section{Proofs of the main theorems}\label{sec:deduction}

We first establish the following quantitative estimate, allowing an auxiliary function $\phi(A)$.

\begin{proposition}\label{prop:exceptional}
Let $\phi:\RR_{>0}\to\RR_{>1}$ satisfy
$\phi(A)\leq(\log A)^{1/2}$ for all sufficiently large $A$. Then
\[
 \frac{\#\{V\in\VV^{\loc}(A):N_V(A\phi(A))=0\}}
      {\#\VV(A)}
 \ll\phi(A)^{-1/72},
\]
with an implied constant independent of $\phi$.
\end{proposition}

\begin{proof}
Put $B=A\phi(A)$. By hypothesis, $A<B\leq A(\log A)^{1/2}$ for all sufficiently large $A$, so $A$ and $B$ lie in the range~\eqref{eq:range}. Since $\#\{\a\in\ZZ^{20}:\|\a\|<A^{7/8}\}\ll A^{35/2}$ and $\#\mathbb{V}(A)\asymp A^{20}$, the surfaces $V\in\mathbb{V}(A)$ with $\|\a_V\|<A^{7/8}$ have a proportion which is $O(A^{-5/2})$, and we discard them. For the remaining surfaces,~\eqref{eq:local-loss} gives
\[
 0\leq N_V^{\loc}(B)-L_{\a_V}^{\good}(B,A^{5/12})=o(1).
\]
If $N_V(B)=0$, then $N_{\a_V}^{\good}(B,A^{5/12})=0$ as well.
Thus, for all sufficiently large $A$, such a surface satisfies either
\[
 N_V^{\loc}(B)\leq\phi(A)^{2/3}
\]
or
\[
 |N_{\a_V}^{\good}(B,A^{5/12})-L_{\a_V}^{\good}(B,A^{5/12})|
                          >\tfrac12\phi(A)^{2/3}.
\]
For locally soluble surfaces, the first alternative has proportion
$O(\phi(A)^{-1/72})$ by~\eqref{eq:local-small}. For the second, apply
Chebyshev's inequality and Proposition~\ref{prop:variance} to bound this 
proportion by
\[
 \ll\phi(A)^{-4/3}\cdot\phi(A)=\phi(A)^{-1/3}.
\]
Consequently,
\[
 \frac{\#\{V\in\VV^{\loc}(A):N_V(A\phi(A))=0\}}
      {\#\VV(A)}
 \ll A^{-5/2}+\phi(A)^{-1/3}+\phi(A)^{-1/72}
 \ll\phi(A)^{-1/72},
\]
as required.
\end{proof}

\begin{proof}[Proof of Theorem~\ref{thm:main}]
Let 
$\phi(A)=(\log A)^{1/2}$.
Every locally soluble surface with no
rational point is counted 
by the $V\in\VV^{\loc}(A)$ for which $N_V(A\phi(A))=0$.
Applying Proposition~\ref{prop:exceptional}, 
we see that  the proportion of
surfaces failing the Hasse principle is $O((\log A)^{-1/144})$.
Since rational solubility implies local solubility, it follows that
\[
 0\leq\frac{\#\VV^{\loc}(A)}{\#\VV(A)}-\rho(A)
 \ll(\log A)^{-1/144}.
\]
The first term tends to $\rho^{\loc}$
by~\cite[Theorem~3.6]{poonen}, proving the theorem.
\end{proof}

\begin{proof}[Proof of Theorem~\ref{thm:height}]
Fix $0<\eta<1$. The hypothesis on $\psi$ implies that, for all
sufficiently large $A$ and every $V$ with
$\eta A<\norm{\a_V}\leq A$, we have
\[
 \psi(\norm{\a_V})\geq\frac{\norm{\a_V}}{\eta^2}
 >\frac A\eta.
\]
For such a surface, $\mathfrak M(V)>\psi(\norm{\a_V})$ implies
$N_V(A/\eta)=0$. Proposition~\ref{prop:exceptional}, with the
constant choice $\phi=1/\eta$, therefore bounds the proportion of
locally soluble surfaces in this range which violate the desired
height bound by $O(\eta^{1/72})$. This choice of $\phi$ is admissible
once $A$ is sufficiently large in terms of $\eta$.
Surfaces of coefficient height at most $\eta A$ account for a
proportion $O(\eta^{20})$. Thus
\[
 \limsup_{A\to\infty}
 \frac{\#\{V\in\VV^{\loc}(A):
             \mathfrak M(V)>\psi(\norm{\a_V})\}}
      {\#\VV(A)}
 \ll\eta^{1/72}+\eta^{20}.
\]
Letting $\eta\to0$ proves that this proportion tends to zero.
Every surface satisfying the height bound is locally soluble, so
the asserted limit follows from the existence of
$\rho^{\loc}$.
\end{proof}

\section{Brauer--Manin obstructions}\label{s:BMO}

In this section we prove Theorem~\ref{thm:BM}. 
For 
a smooth cubic surface $V\subset \PP^3$ over $\QQ$,
with structure
morphism $\pi\colon V\to\operatorname{Spec}\QQ$, we write $\Br_0(V)=\pi^*\Br(\QQ)$ for the image of
$\Br(\QQ)\to\Br(V)$. Theorem~\ref{thm:BM} is based on  the following result, where  $\mathbb V^{\mathrm{sm}}(A)$ denotes the set of smooth
surfaces in $\mathbb V(A)$.

\begin{proposition}\label{prop:Br}
We have
\[
\#\bigl\{V\in\mathbb{V}^{\mathrm{sm}}(A):\ \Br(V)/\Br_0(V)\neq 0\bigr\}\ll A^{19}.
\]
\end{proposition}

\begin{proof}[Proof of Theorem~\ref{thm:BM} assuming Proposition~\ref{prop:Br}]
Let $V\in\mathbb{V}^{\mathrm{loc}}(A)$ be smooth such that $\Br(V)/\Br_0(V)=0$. Then any
$\gamma\in\Br(V)$ is of the form $\gamma=\pi^*\beta$, with $\beta\in\Br(\QQ)$. Thus 
\[
\sum_v \operatorname{inv}_v\gamma(x_v)=\sum_v\operatorname{inv}_v\beta=0,
\]
for any $(x_v)_v\in V(\mathbb{A}_{\QQ})$, 
by the reciprocity law of global class field theory
\cite{Tate}. Thus constant classes impose no Brauer--Manin obstruction and we have 
$V(\mathbb{A}_{\QQ})^{\Br}=V(\mathbb{A}_{\QQ})\neq\emptyset$. It follows that every smooth
surface counted in Theorem~\ref{thm:BM} satisfies $\Br(V)/\Br_0(V)\neq0$, and there are
$O(A^{19})$ of these by Proposition~\ref{prop:Br}.

The singular surfaces in $\mathbb{V}(A)$ correspond to coefficient vectors
$\mathbf{a}\in\ZZp^{20}$ with $\norm{\mathbf{a}}\le A$ and $\Disc(\mathbf{a})=0$, where
$\Disc\in\ZZ[\mathbf{a}]$ denotes the discriminant of a cubic form in four variables.
But a fixed nonzero polynomial in twenty variables has $O(A^{19})$ integral zeros in a box of side $A$.
The theorem therefore follows from the observation that 
$\#\mathbb V^{\mathrm{loc}}(A)\asymp A^{20}$, by \cite[Theorem~3.6]{poonen}.
\end{proof}

We expect the exponent $19$ in Proposition~\ref{prop:Br} to be far from the
truth. Indeed, a non-constant Brauer group forces the Galois group of the
$27$ lines to be a proper subgroup of $W(E_6)$, and O'Dorney \cite[Conjecture~1.2]{ODorney}
conjectures that this happens for only $O_{\varepsilon}(A^{16+\varepsilon})$ surfaces in $\VV(A)$.

The rest of this section is devoted to the proof of Proposition~\ref{prop:Br}. We first
recall why a non-constant Brauer group forces the Galois action on the $27$ lines to fix a
double-six or a Steiner decomposition. This is summarised in Lemma~\ref{lem:EJ}.
We then construct covers
\[
\mathcal D,\mathcal T\longrightarrow U
\]
of the space \(U\) of smooth cubic surfaces such that the parameter
point of every such surface lifts to \(\mathcal D(\mathbb Q)\) or
\(\mathcal T(\mathbb Q)\). We show that these covers are integral of
degrees \(36\) and \(40\), respectively.  We finally conclude with an application of work by  Buggenhout, Cluckers, Salberger, Santens and Vermeulen on thin sets of type~II \cite{BCSV}.

\subsection*{Brauer groups of cubic surfaces}

Throughout this subsection, $S$ denotes a smooth cubic surface over $\QQ$;
we write $\Sbar=S\times_{\QQ}\overline{\QQ}$ and
$\Gamma_{\QQ}=\Gal(\overline{\QQ}/\QQ)$.

\subsubsection*{Lines and markings}
The $27$ lines on $\Sbar$ have the same incidence graph for every smooth
cubic surface. Fix an abstract model $\mathcal G$ of this graph, with
vertex set $\Lzero$, where adjacency corresponds to intersection of
distinct lines. We identify
$\operatorname{Aut}(\mathcal G)$ with the Weyl group $W(E_6)$, which has
order $51840$ by  \cite[Chapter~4; Remark~1.7 and Section~3.6]{Huybrechts}.

A \emph{marking} of $S$ is an incidence-preserving bijection $\lambda$
from $\Lzero$ to the lines on $\Sbar$. The group $W(E_6)$ acts simply
transitively on markings by
$w\cdot\lambda=\lambda\circ w^{-1}$.
Galois acts by $\sigma\cdot\lambda=\sigma\circ\lambda$,
and these two actions commute. Here Galois acts on the lines on
$\Sbar$, while the reference labels in $\Lzero$ remain fixed.
A marking identifies the Galois action on the lines with the homomorphism
$$
\rho_\lambda\colon\Gamma_{\QQ}\longrightarrow W(E_6),
\qquad
\rho_\lambda(\sigma)=\lambda^{-1}\circ\sigma\circ\lambda.
$$
Replacing $\lambda$ by $\lambda\circ w$ replaces $\rho_\lambda$
by $w^{-1}\rho_\lambda w$. Thus the image
$G_S=\rho_\lambda(\Gamma_{\QQ})$ is well defined up to conjugacy
in $W(E_6)$.

\subsubsection*{The algebraic Brauer group}
As explained in 
\cite[Lemma~2.5]{EJ-order3}, for example, there is a canonical isomorphism 
$\Br(S)/\Br_0(S)\cong H^1(\Gamma_{\QQ},\Pic(\Sbar)).
$ A theorem of Swinnerton-Dyer \cite{SD93} shows that this group is one of
\[
0,\qquad\ZZ/2\ZZ,\qquad(\ZZ/2\ZZ)^2,\qquad\ZZ/3\ZZ,\qquad(\ZZ/3\ZZ)^2.
\]
In particular, only $2$-torsion and $3$-torsion occur.

\subsubsection*{Two configurations of lines}
A \emph{double-six} is an unordered pair
$$
\{\{a_1,\dots,a_6\},\{b_1,\dots,b_6\}\}$$ 
of disjoint sets of six
pairwise skew lines, labelled so that $a_i$ meets $b_j$ if and only
if $i\neq j$.
A \emph{pair of Steiner trihedra} may be identified with a set of
nine distinct lines that can be arranged in a $3\times3$ array
so that each row and each column forms a tritangent plane section;
see \cite[Section 3]{EJ-order3}.
The choice of array is not part of the configuration.
A \emph{Steiner decomposition} is an unordered partition of the
$27$ lines into three such sets of nine lines.
Three distinct lines on a smooth cubic surface lie in a tritangent
plane if and only if they pairwise meet by
\cite[Chapter~4, Section~3.3]{Huybrechts}.
Thus both configurations are determined by incidence and are
preserved by graph isomorphisms.

We write $\mathcal D(S)$ and $\mathcal T(S)$ for the sets of
double-sixes and Steiner decompositions on $\Sbar$.
Likewise, write $\mathcal D(\Lzero)$ and $\mathcal T(\Lzero)$
for the sets of these configurations defined using adjacency
in the fixed graph $\mathcal G$.
These sets have cardinalities
$36$ and $40$, respectively, and $W(E_6)$ acts transitively on
each of the two model sets, by 
\cite[\S3.1]{EJ-BM} and \cite[Remark~3.4(ii)]{EJ-order3}.
Fix $D_0\in\mathcal D(\Lzero)$ and $T_0\in\mathcal T(\Lzero)$,
and let $H_D$ and $H_T$ be their stabilisers. Then
\[
[W(E_6):H_D]=36,\qquad [W(E_6):H_T]=40.
\]

A configuration is \emph{$\Gamma_{\QQ}$-invariant} if it is
preserved as a whole. Thus Galois may interchange the two halves
of an invariant double-six, or permute the three parts of an
invariant Steiner decomposition.
A marking $\lambda$ transports configurations from the model
to $\Sbar$, with
$\sigma(\lambda(C))
 =\lambda\bigl(\rho_\lambda(\sigma)C\bigr)$.
By transitivity, this gives
\begin{equation}\label{eq:inv-conj}
\mathcal D(S)^{\Gamma_{\QQ}}\neq\emptyset
\iff
G_S\subseteq wH_Dw^{-1}
\quad\text{for some }w\in W(E_6),
\end{equation}
and similarly for $\mathcal T(S)$ and $H_T$.

The following consequence of the results of Elsenhans and Jahnel
is the Brauer-group input to our counting argument.

\begin{lemma}\label{lem:EJ}
Let $S$ be a smooth cubic surface over $\QQ$ with
$\Br(S)/\Br_0(S)\neq0$. Then
\[
\mathcal D(S)^{\Gamma_{\QQ}}\neq\emptyset
\quad\text{or}\quad
\mathcal T(S)^{\Gamma_{\QQ}}\neq\emptyset.
\]
Equivalently, $G_S$ is contained in a conjugate of $H_D$ or $H_T$.
\end{lemma}

\begin{proof}
By the identification with
$H^1(\Gamma_{\QQ},\Pic(\Sbar))$ and the classification recalled
above, $\Br(S)/\Br_0(S)$ has order $2$, $4$, $3$ or $9$.
In the first two cases,
\cite[Theorem~5.11(b)(i)]{EJ-BM}
gives a $\Gamma_{\QQ}$-invariant double-six.
In the latter two cases,
\cite[Proposition~3.15]{EJ-order3}
gives a Galois-invariant ordered triplet of pairs of Steiner
trihedra. Forgetting its ordering gives a
$\Gamma_{\QQ}$-invariant Steiner decomposition.
The final assertion follows from \eqref{eq:inv-conj}
and its analogue for $\mathcal T(S)$.
\end{proof}

\subsection*{The configuration covers}

Let $U\subset\PP^{19}_{\QQ}$ be the open subset parametrising
smooth cubic surfaces. For a field $K\supseteq\QQ$ and
$u\in U(K)$, write $S_u$ for the corresponding cubic surface.
For a morphism $X\to U$, write $X_u$ for its fibre above $u$.

\begin{proposition}\label{prop:integral}
There are integral schemes $\mathcal D$ and $\mathcal T$ with
finite étale morphisms
$$
p_{\mathcal D}\colon\mathcal D\to U,
\qquad
p_{\mathcal T}\colon\mathcal T\to U
$$
of degrees $36$ and $40$, respectively.
Any $u\in U(\QQ)$ lifts to $\mathcal D(\QQ)$, respectively
$\mathcal T(\QQ)$, if and only if $S_u$ has a
$\Gamma_{\QQ}$-invariant double-six, respectively Steiner
decomposition.
\end{proposition}

\begin{proof}
Let $\mathcal L\to U$ be the scheme parametrising pairs $(L,u)$,
with $L$ a line on $S_u$. This morphism is projective and smooth
of relative dimension zero by
\cite[Chapter~2, Corollary~1.14]{Huybrechts}.
It is therefore finite étale, of degree $27$.

Let $U'\to U$ be a connected finite \'etale cover over which $\mathcal L$
splits into $27$ copies of $U'$, so that the $27$ lines can be labelled.
For instance, one can take $U'$ to be a connected component of the complement
of all the diagonals in the $27$-fold fibre product of $\mathcal L$ over $U$.
We show that the incidence matrix of the labelled lines is constant on $U'$.
Since $U'$ is smooth and connected, it is integral.
Any two labelled lines that meet at the generic point of $U'$
meet in every fibre, since meeting is a closed condition.
There are exactly $270$ ordered pairs
of distinct meeting lines in every geometric fibre: each of
the $27$ lines meets ten others by
\cite[Chapter~4, Section~3.3]{Huybrechts}.
Thus no additional pair can meet in a special fibre, proving
the assertion.

Let $\mathcal L^{\Lzero}_U$ be the fibre product of $27$ copies of $\mathcal L$ over $U$, indexed
by $\Lzero$, so that for an algebraically closed field $K$ and $u\in U(K)$ a point of
$(\mathcal L^{\Lzero}_U)_u(K)$ is a family $(L_\ell)_{\ell\in\Lzero}$ of lines on $S_u$, not
necessarily distinct. 
Let $\mathcal M$
be its locus of families of distinct lines whose incidence
graph agrees with $\mathcal G$ under the given labelling.
To see that $\mathcal M$ is open and closed, we base change to the cover $U'\to U$ introduced above. The product
$\mathcal L^{\Lzero}_U\times_U U'$ is a disjoint union
of copies of $U'$, each corresponding to a fixed choice
of a line for every label $\ell\in\Lzero$.
On each copy, distinctness is constant, and the preceding
paragraph shows that incidence is constant.
Thus the conditions defining $\mathcal M$ select a union
of these copies.
Since $U'\to U$ is finite \'etale and surjective, openness and closedness can be checked after this base change. Hence $\mathcal M$ is open and closed in $\mathcal L^{\Lzero}_U$.
It is therefore finite étale over $U$, and its geometric
fibres are precisely the sets of markings.
Permuting the factors gives an action of $W(E_6)$ over $U$
that preserves $\mathcal M$. On markings this is
$w\cdot\lambda=\lambda\circ w^{-1}$, so the action is
simply transitive on each geometric fibre.

For a subgroup $H\subseteq W(E_6)$, the quotient $\mathcal M/H$ is finite
\'etale over $U$. Indeed, after labelling the lines it is a disjoint
union of copies of the base, one for each $H$-orbit of markings.
Each orbit contains $|H|$ markings, so the degree is $[W(E_6):H]$.
Set
$$
\mathcal D=\mathcal M/H_D,\qquad
\mathcal T=\mathcal M/H_T.
$$
Fix $u\in U(\QQ)$. A marking $\lambda$ sends the reference labels
in $\Lzero$ to the lines on $\overline{S_u}$. The map
$$
 \mathcal D_u(\overline{\QQ})\longrightarrow\mathcal D(S_u),
 \qquad [\lambda]\longmapsto\lambda(D_0)
$$
is a bijection. Every double-six occurs by transitivity, and two
markings give the same double-six precisely when they differ by an
element of $H_D$, the stabiliser of $D_0$. The map commutes with
Galois, so a rational lift of $u$ is the same as a Galois-invariant
double-six. The same argument applies to $T_0$ and $H_T$.
(By comparison, a rational point of $\mathcal M$ above $u$
amounts to a Galois-invariant marking, which requires all
$27$ lines to be rational.)

It remains to prove integrality. The geometric monodromy
group of the $27$ lines is $W(E_6)$ by
\cite[Chapter~4, Corollary~1.14]{Huybrechts}.
The Galois group over $\QQ(U)$ contains this group and is
itself contained in $W(E_6)$, so it is also $W(E_6)$.
It therefore acts transitively on the markings of the generic
surface, and hence on their $H_D$-orbits and $H_T$-orbits.
It follows that the generic fibres of $\mathcal D$ and
$\mathcal T$ are connected. Indeed, the connected components
of a finite étale scheme over a field correspond to the
Galois orbits on its geometric points.
Both covers are normal, since they are finite étale over
the smooth scheme $U$, and every connected component
dominates $U$. Their connected generic fibres therefore imply
that they are connected and hence integral.
\end{proof}

\subsection*{Type~II thin sets}

Let $f\colon X\to\PP^n$ be a morphism of varieties over $\QQ$. Recall that $f$ is
\emph{quasi-finite} if all its fibres are finite, and \emph{dominant} if its image is Zariski
dense. 
If $f$ is quasi-finite and dominant and $X$ is integral,
then $\dim X=n$, and the \emph{degree} of $f$ is
the degree $[\QQ(X):\QQ(\PP^n)]$ of the extension of function fields. We use the following result.

\begin{theorem}[{\cite[Theorem~1.1]{BCSV}}]\label{thm:BCSV}
Let $f\colon X\to\PP^n$ be a quasi-finite dominant morphism over $\QQ$ of degree $d\ge5$, with
$X$ integral. Then, for all $B\ge2$,
\[
\#\{x\in f(X(\QQ)) : H(x)\le B\}\ll_f B^{n},
\]
where $H$ is the standard exponential height on $\PP^n(\QQ)$.
\end{theorem}

\begin{proof}[Proof of Proposition~\ref{prop:Br}]
Let $V\in\mathbb V^{\mathrm{sm}}(A)$ with
$\Br(V)/\Br_0(V)\neq0$, and let $u=[\mathbf{u}]\in U(\QQ)$ be the
corresponding point, so that $V=S_u$ and $H(u)\leq\norm{\mathbf u}\le A$, where $H$ is the height
in Theorem~\ref{thm:BCSV}.
By Lemma~\ref{lem:EJ}, the surface $S_u$ has a
$\Gamma_{\QQ}$-invariant double-six or Steiner decomposition.
Proposition~\ref{prop:integral} therefore gives
$$
u\in
p_{\mathcal D}\bigl(\mathcal D(\QQ)\bigr)
\cup
p_{\mathcal T}\bigl(\mathcal T(\QQ)\bigr).
$$
By the same proposition, $\mathcal D$ and $\mathcal T$ are
integral. Moreover, on
composing these finite surjective morphisms with
the open immersion $U\hookrightarrow\PP^{19}_{\QQ}$, we deduce that 
their maps to $\PP^{19}_{\QQ}$ are quasi-finite
and dominant, of degrees $36$ and $40$, respectively.
Applying Theorem~\ref{thm:BCSV} to each map gives
$O(A^{19})$ points of height at most $A$ in this union,
proving the proposition.
\end{proof}


\begin{thebibliography}{99}

\bibitem{chris}
L. Beneish and C. Keyes,
How often does a cubic hypersurface have a rational point? {\em Selecta Math.} {\bf 31}:92 (2025).

\bibitem{Bhargava}
M.~Bhargava,
A positive proportion of plane cubics fail the Hasse principle.
{\em Preprint}, 2014. 
\newline
(\href{https://arxiv.org/abs/1402.1131}{arXiv:1402.1131})

\bibitem{BLS}
T.~Browning, P.~Le Boudec and W.~Sawin,
The Hasse principle for random Fano hypersurfaces.
{\em Annals of Math.} {\bf 197} (2023), 1115--1203.

\bibitem{BCSV} T.~Buggenhout, R.~Cluckers, P.~Salberger, T.~Santens and F.~Vermeulen,
Serre's question on thin sets in projective space. {\em Preprint}, 2025. 
(\href{https://arxiv.org/abs/2506.13471}{arXiv:2506.13471})

\bibitem{CTS} 
J.-L. Colliot-Th\'el\`ene and J.-J. Sansuc,
La descente sur les vari\'et\'es rationnelles. In {\em Journ\'ees de 
g\'eom\'etrie alg\'ebrique d'Angers (Juillet 1979)}, 223--237, 
Sijthoff \& Noordhoff, 1980. 


\bibitem{EJ-BM} A.-S.~Elsenhans and J.~Jahnel, On the Brauer--Manin obstruction
for cubic surfaces. \emph{J.\ Comb.\ Number Theory} \textbf{2} (2010), 107--128.

\bibitem{EJ-order3} A.-S.~Elsenhans and J.~Jahnel, On the order three Brauer
classes for cubic surfaces. \emph{Cent.\ Eur.\ J.\ Math.} \textbf{10} (2012), 903--926.

\bibitem{Huybrechts} D.~Huybrechts, \emph{The geometry of cubic hypersurfaces}. Cambridge
Studies in Advanced Mathematics \textbf{206}, Cambridge University Press, Cambridge, 2023.

\bibitem{LeBoudec}
P.~Le Boudec,
Height of rational points on random Fano hypersurfaces.
{\em Algebra \& Number Theory} {\bf 15} (2021), 657--672.


\bibitem{ODorney}
E. O'Dorney,
Van der Waerden's problem for cubic surfaces.
{\em Preprint}, 2026. (\href{https://arxiv.org/abs/2610.37945}{arXiv:2610.37945})

\bibitem{poonen}
B.~Poonen and J.~F.~Voloch,
Random Diophantine equations. In
{\em Arithmetic of higher-dimensional algebraic varieties
(Palo Alto, CA, 2002)}, 175--184,
Progr. Math. {\bf 226}, Birkh\"auser Boston, Boston, MA, 2004.

\bibitem{SD93} P.~Swinnerton-Dyer, The Brauer group of cubic surfaces. \emph{Math.\ Proc.\
Cambridge Philos.\ Soc.} \textbf{113} (1993), 449--460.

\bibitem{Tate} J.~T.~Tate, Global class field theory.  In  \emph{Algebraic Number Theory
(Proc. Instructional Conf., Brighton, 1965)}, 162--203, J.~W.~S.~Cassels and A.~Fr\"ohlich (eds.),
Academic Press, London; Thompson Book Co., Washington, D.C., 1967.

\bibitem{anna}
A. Zanoli,
Average density of rational points on quartic surfaces.
{\em Preprint}, 2026.
\newline
(\href{https://arxiv.org/abs/2611.32640}{arXiv:2611.32640})

\end{thebibliography}
\end{document}